\documentclass[11pt]{article}

\usepackage{amssymb}
\usepackage{amsmath}
\usepackage{shadow}
\usepackage{fancybox}

\def\0{\mbox{\tiny $0$}}
\def\1{\mbox{\tiny $1$}}
\def\2{\mbox{\tiny $2$}}
\def\3{\mbox{\tiny $3$}}
\def\4{\mbox{\tiny $4$}}
\def\5{\mbox{\tiny $5$}}
\def\6{\mbox{\tiny $6$}}
\def\7{\mbox{\tiny $7$}}
\def\8{\mbox{\tiny $8$}}
\def\9{\mbox{\tiny $9$}}

\def\R{\mathbb{R}}

\def\H{\mathbb{H}}

\def\ij2{\mbox{\tiny $\frac{i+j}{2}$}}
\def\ji2{\mbox{\tiny $\frac{j-i}{2}$}}

\begin{document}

\title{{\huge {\bf Comments on the Matrix Representations of
Quaternions}}}
\author{
{\Large \sf  Gisele Ducati}\\
{\em Departamento de Matem\'atica, UFPR}\\
{\small  CP 19081, 81531-970  Curitiba (PR) Brasil}\\
{\small  \tt ducati@mat.ufpr.br} }

\maketitle










\begin{quote}
{\bf Abstract} ~Some comments are made on the matrices which serve
as the basis of a quaternionic algebra. We show that these
matrices are related with the quaternionic action of the imaginary
units from the left and from the right.

\noindent {\bf MSC.} 11R52 - 15A33.
\end{quote}


\maketitle




We start by recalling the notation of the referred
paper~\cite{FGT03}. It's well known that a quaternion
\[
q = a + b h + c j + d k~,~~~~a,b,c,d \in \mathbb{R}~,
\]
where $h, j$ and $k$ denote the imaginary units of the algebra of
quaternions
\[
h^2 = j^2 = k^2 = h j k = -1~,
\]
can be written in matrix form as
\[
Q = a I + b H + c J + d K~,
\]
with
\[
H^2 = J^2 = K^2 = H J K = - I~.
\]
Like pointed out in the paper~\cite{FGT03} the matrices $H, J$ and
$K$ can be chosen in different ways so adopting skew-symmetric
signed permutation matrices with a plus one in its first row we
find six signed permutation matrices which arrange themselves to
form two Hamiltonian systems, given by $H_1, J_1, K_1$ and $H_2,
J_2, K_2$.

We emphasize that the first system represents the left action of
the imaginary units $j, h$ and $k$, respectively, and the matrices
of the second system represents, respectively, the right action of
the quaternionic imaginary units, $h, j$ and $k$. In both cases we
should multiply the units by $- 1$. Let $q = q_0 + h q_1 + j q_2 +
k q_3$ and $p = p_0 + h p_1 + j p_2 + k p_3$ real quaternions. If
we distinguish between the left and right action of the
quaternionic imaginary units $h$, $j$, and $k$ by introducing the
operators $L_{q}$ and $R_{p}$ whose action on quaternionic
functions $\Psi$, $\Psi: \H \to \R$, is given by
\begin{equation}
L_{q} \Psi = q \, \Psi~~~\mbox{and}~~~ R_{p} \Psi =\Psi \, p~,
\end{equation}
where
\[
L_h^2 = L_j^2 = L_k^2 = R_h^2 = R_j^2 = R_k^2 = L_h L_j L_k = R_k
R_j R_h = -1~,
\]
and
\begin{equation}
\left[ \, L_{q} \, , \, R_{p} \, \right] = L_q\,R_p - R_p \, L_q =
0~,
\end{equation}
we find the following {\em real} matrix
representation~\cite{DD99,DD01}
\begin{equation}
\label{rrep} L_{q} \leftrightarrow
 \left( \begin{array}{rrrr} q_{\0} & $-$
q_{\1} & $-$
q_{\2} & $-$ q_{\3} \\
q_{\1} & q_{\0} & $-$
q_{\3} & q_{\2} \\
q_{\2} & q_{\3} &
q_{\0} & $-$ q_{\1} \\
q_{\3} & $-$ q_{\2} & q_{\1} &  q_{\0} \end{array}
\right)~~~\mbox{and}~~~R_{p} \leftrightarrow  \left(
\begin{array}{rrrr} p_{\0} & $-$ p_{\1} & $-$
p_{\2} & $-$ p_{\3} \\
p_{\1} & p_{\0} &
p_{\3} & $-$ p_{\2} \\
p_{\2} & $-$ p_{\3} &
p_{\0} &  p_{\1} \\
p_{\3} &  p_{\2} & $-$ p_{\1} &  p_{\0} \end{array} \right)~.
\end{equation}
It's interesting to note that the six matrices
\begin{equation}
\label{op} L_h,~ L_j,~ L_k,~ R_h,~ R_j~~\mbox{and}~~R_k
\end{equation}
are the generators of the one-dimensional quaternion unitary group
linear from the right. See \cite{DD99,DUC02} for further details.
In addition, it's more natural to arrange the six signed
permutation matrices cited in \cite{FGT03} as given in (\ref{op}).

Now, let us see how many triplets, obeying the quaternionic
algebra, is possible to write with these elements. We should
consider each element from (\ref{op}) and it's conjugate (or it's
negation) which give 12 elements. Obviously, fixing a imaginary
unit, the second one must be chosen in a way that it does not
commute with the previously fixed one, otherwise, the
multiplication of the imaginary units will not square to minus
one, contradicting the quaternion algebra. So, if $L_h$ has been
chosen for the first imaginary unit, the remaining possibilities
are $\pm L_j$ and $\pm L_k$ and the third imaginary unit is
determined by the multiplication of the first two. So we have $12$
possibilities for the first imaginary unit and $4$ for the second,
that means, $12 \times 4 \times 1$, totalizing 48 possibilities.
Of course, we can consider all automorphisms of the quaternion
algebra, or alternatively, thinking in a geometric way, any
quaternion can be rotated in the three dimensional space, then we
should consider the transformation
\[
q \to q_u = u q u^{-\1}~,~~~|u|=1~.
\]
Note that
\[
h_u^2 = j_u^2 = k_u^2 = h_u j_u k_u = -1~.
\]
This similarity transformation can be made with $L_h,~ L_j,~ L_k,~
R_h,~ R_j$ and $R_k$. The automorphism group of the units is the
$O(3)$ subgroup of $O(4)$ given by
\[
\Omega' = U \Omega U^{\dag}~.
\]

\end{document}